\documentclass[11pt]{article}
\usepackage{amssymb, amsmath, amsthm, mathrsfs, ascmac, color}
\usepackage[pdftex]{graphicx} 
\usepackage[subrefformat=parens]{subcaption}
\usepackage[top=3cm, bottom=3cm, left=3cm,right=3cm]{geometry}
\usepackage{comment}

\theoremstyle{plain}
\newtheorem{thm}{Theorem}[section]

\newtheorem{rmk}{Remark}
\makeatletter

\@addtoreset{equation}{section}
\makeatother

\newcommand{\dd}{\mathrm d}
\newcommand{\pd}{\partial}
\newcommand{\ee}{\mathrm e}
\newcommand{\EE}{\mathbb E}
\newcommand{\PP}{\mathbf{P}}
\newcommand{\OO}{\mathcal O}

\newcommand{\cov}{\mathrm{Cov}}

\newcommand{\pto}{\stackrel{\PP}{\to}}

\allowdisplaybreaks
\title{\textbf{Estimation for the damping factor 
of the driving process of 
an SPDE in two space dimensions 
}}
\date{}

\author{\textbf{Yozo Tonaki}\thanks{Graduate School of Engineering Science, Osaka University}
\thanks{Center for Mathematical Modeling and Data Science (MMDS), Osaka University} 
\footnote{e-mail: \texttt{y.tonaki.es@osaka-u.ac.jp}}
\and \textbf{Yusuke Kaino}\thanks{Graduate School of Maritime Sciences, Kobe University}
\and \textbf{Masayuki Uchida}$^{* \dag}$\thanks{CREST, Japan Science and Technology Agency}
}

\begin{document}
\bibliographystyle{plain}
\maketitle

\begin{abstract}
We study parametric estimation 
for a second order linear parabolic stochastic partial differential equation (SPDE)
in two space dimensions driven by a $Q$-Wiener process
based on high frequency spatio-temporal data.
We give an estimator of the damping parameter of the $Q$-Wiener process  
of the SPDE based on quadratic variations with 
temporal and spatial increments.
We also provide simulation results of the proposed estimator.

\begin{center}
\textbf{Keywords and phrases}
\end{center}
Damping factor,
high frequency spatio-temporal data,
linear parabolic stochastic partial differential equations,
parametric estimation, 
$Q$-Wiener process.
\end{abstract}

\section{Introduction}
We consider the following linear parabolic stochastic partial differential equation
(SPDE) in two space dimensions
\begin{align}
\dd X_t(y,z)
&=\biggl\{
\theta_2
\biggl(\frac{\pd^2}{\pd y^2}+\frac{\pd^2}{\pd z^2}\biggr)
+\theta_1\frac{\pd}{\pd y} 
+\eta_1\frac{\pd}{\pd z} 
+\theta_0 
\biggr\} X_t(y,z) \dd t
\nonumber
\\
&\qquad
+\sigma\dd W_t^{Q}(y,z),
\quad (t,y,z) \in [0,1] \times D
\label{2d_spde}
\end{align}
with the initial condition $X_0 = \xi$ and the Dirichlet boundary condition
\begin{equation*}
X_t(y,z) = 0, \quad (t,y,z) \in [0,1] \times \pd D,
\end{equation*}
where $D =(0,1)^2$, $W_t^Q$ is a $Q$-Wiener process in a Sobolev space on $D$,
the initial value $\xi$ is independent of $W_t^Q$,
and $(\theta_0, \theta_1, \eta_1, \theta_2, \sigma) \in \mathbb R^3 \times (0,\infty)^2$ 
are unknown parameters.

Statistical inference for second order parabolic SPDEs based on high frequency data
has been studied by some researchers, see for example, 
Bibinger and Trabs \cite{Bibinger_Trabs2020},
Chong \cite{Chong2020},
Cialenco and Huang \cite{Cialenco_Huang2020},
Kaino and Uchida \cite{Kaino_Uchida2020,Kaino_Uchida2021},
Hildebrandt and Trabs \cite{Hildebrandt_Trabs2021},
Bibinger and Bossert \cite{Bibinger_Bossert2023},
Tonaki et al.\,\cite{TKU2023,TKU2024a,TKU2024b,TKU2024arXiv},
and Gaudlitz and Reiss \cite{Gaudlitz_Reiss2023}.
Recently, Bossert \cite{Bossert2023arXiv} considered parametric estimation for 
second order parabolic SPDEs in multiple space dimensions.
In particular, \cite{Bossert2023arXiv} constructed an estimator of 
the damping parameter $\alpha$ of a $Q$-Wiener process based on realized volatilities,
and stated that the estimator of the damping parameter 
is asymptotically normal with the rate 
$1/\sqrt{mn}$, where $n$ and $m$ denote the number of temporal and spatial observations, 
respectively, and $m = \OO(n^\rho)$ with for some $\rho \in (0,(1-\alpha)/4)$
($\alpha \in (0,1)$) in the case of the SPDE in two space dimensions.

In this paper, we propose an estimator of the damping parameter $\alpha$ 
of the $Q$-Wiener process \eqref{QW} below 
based on temporal and spatial quadratic variations.
The main purpose of this paper is that the estimator of $\alpha$ is 
bounded in probability at the parametric rate $1/\sqrt{m N}$,
where $m$ is the number of spatially thinned data, 
$N$ is the number of temporal observations, $m = \OO(N)$ and $N = \OO(m)$.

This paper is organized as follows.
In Section \ref{sec2}, we provide an estimator of the damping parameter $\alpha$
of the $Q$-Wiener process in SPDE \eqref{2d_spde}, and show the asymptotic properties 
of the estimator of $\alpha$.
In Section \ref{sec3}, we give the simulation results of the proposed estimator.
Section \ref{sec4} is devoted to the proof of the result of Section \ref{sec2}.

\section{Main results}\label{sec2}
\subsection{Setting and notation}
We define an operator $A_\theta$ by 
\begin{equation*}
-A_\theta = 
\theta_2\biggl(\frac{\pd^2}{\pd y^2} + \frac{\pd^2}{\pd z^2} \biggr)
+ \theta_1\frac{\pd}{\pd y} + \eta_1\frac{\pd}{\pd z} + \theta_0.
\end{equation*}
The eigenfunctions $e_{l_1,l_2}$ of the operator $A_\theta$ 
and the corresponding eigenvalues $\lambda_{l_1,l_2}$ are given by 
\begin{equation*}
e_{l_1,l_2}(y,z)=e_{l_1}^{(1)}(y)e_{l_2}^{(2)}(z),
\quad
\lambda_{l_1,l_2}=\theta_2(\pi^2(l_1^2+l_2^2)+\Gamma)
\end{equation*}
for $l_1,l_2 \ge 1$ and $y,z\in [0,1]$, where
\begin{equation*}
e_{l_1}^{(1)}(y)
=\sqrt 2 \sin(\pi l_1 y) \ee^{-\kappa y/2},
\quad
e_{l_2}^{(2)}(z)
=\sqrt 2 \sin(\pi l_2 z) \ee^{-\eta z/2},
\end{equation*}
\begin{equation*}
\kappa=\frac{\theta_1}{\theta_2}, 
\quad
\eta=\frac{\eta_1}{\theta_2},
\quad
\Gamma=-\frac{\theta_0}{\theta_2} +\frac{\kappa^2+\eta^2}{4}.
\end{equation*}
We assume $\lambda_{1,1} > 0$ so that $A_\theta$ is a positive definite and 
self-adjoint operator.
The eigenfunctions $\{e_{l_1,l_2}\}_{l_1,l_2\ge1}$ 
are orthonormal with respect to the weighted $L^2$-inner product
\begin{equation*}
\langle u, v \rangle
= \int_0^1 \int_0^1 u(y,z)v(y,z)\ee^{\kappa y +\eta z} \dd y \dd z, 
\quad 
\| u \| = \sqrt{\langle u, u \rangle}
\end{equation*}
for $u,v \in L^2(D)$.

We consider the $Q$-Wiener process
\begin{equation}\label{QW}
W_t^{Q} = 
\sum_{l_1,l_2\ge1} \mu_{l_1,l_2}^{-\alpha/2} 
e_{l_1,l_2} w_{l_1,l_2}(t)
\end{equation}
with $\mu_{l_1,l_2} = \pi^2(l_1^2 +l_2^2) + \mu_0$, $\mu_0 \in (-2\pi^2,\infty)$ and 
independent real-valued Wiener processes $\{ w_{l_1,l_2} \}_{l_1,l_2 \ge 1}$
for $t \ge 0$ and $\alpha \in (0,2)$, 
where $\mu_0$ and $\alpha$ are unknown parameters,
the parameter space of $\alpha$ is a compact convex subset
of $(0,2)$ and the true value $\alpha^*$ belongs to its interior. 
There exists a unique mild solution of SPDE \eqref{2d_spde}, which is given by 
\begin{equation*}
X_t=\ee^{-t A_\theta}\xi+\sigma\int_0^t \ee^{-(t-s)A_\theta}\dd W_s^{Q},
\end{equation*}
where $\ee^{-t A_\theta} u 
= \sum_{l_1,l_2 \ge 1} \ee^{-\lambda_{l_1,l_2}t}
\langle u, e_{l_1,l_2}\rangle e_{l_1,l_2}$ for $u \in L^2(D)$.
In this paper, we only consider the SPDE driven by the $Q_2$-Wiener process \eqref{QW}, 
but we can discuss the SPDE driven by 
the $Q_1$-Wiener process given in \cite{TKU2023} as well.

We assume that a mild solution is discretely observed on 
the grid $(t_i,y_j,z_k) \in [0,1]^3$ with 
\begin{equation*}
t_i = i\Delta = i/N, \quad 
y_j = j/M_1, \quad
z_k = k/M_2
\end{equation*}
for $i=0,\ldots,N$, $j=0,\ldots,M_1$ and $k=0,\ldots,M_2$.
That is, we have discrete observations 
$\mathbb X_{M,N}=\{X_{t_i}(y_j,z_k)\}$, 
$i=0,\ldots,N$, $j=0,\ldots,M_1$, $k=0,\ldots, M_2$.

For a sequence $\{a_n\}$, 
we will write $a_n \equiv a$ if $a_n = a$ for some $a \in \mathbb R$ and all $n$.

\subsection{Estimation of the damping parameter $\alpha$}
We make the following condition of the initial value $\xi$ 
of SPDE \eqref{2d_spde}.
\begin{description}
\item[\textbf{[A1]}]
The initial value $\xi \in L^2(D)$ is deterministic and
$\| A_\theta^{3/2} \xi \| < \infty$. 
\end{description}

\begin{rmk}
Tonaki et al.\cite{TKU2024b} treated $\alpha \in (0,2)$ as a known parameter 
and assumed $\| A_\theta^{(1+\gamma)/2} \xi \| < \infty$ with $\gamma = \alpha$. 
In this paper, we assume $\| A_\theta^{(1+\gamma)/2} \xi \| < \infty$ 
with $\gamma = 2$ since $\alpha$ is unknown.
\end{rmk}

Fix $b \in (0,1/2)$.
Suppose that we have the following two thinned data of $\mathbb X_{M,N}$.
\begin{align*}
\mathbb X_{m,N}^{(b)} &= \{ X_{t_i}(\widetilde y_j,\widetilde z_k) \},
\quad
i=0,\ldots,N,\ j=0,\ldots,m_1,\ k=0,\ldots, m_2,
\\
\mathbb X_{m',N'}^{(b)} &= \{ X_{t_i'}(\widetilde y_j',\widetilde z_k') \},
\quad
i=0,\ldots,N',\ j=0,\ldots,m_1',\ k=0,\ldots, m_2'
\end{align*}
with $\widetilde y_j = b + j \delta$, $\widetilde z_k = b + k \delta$,
$t_i' = i \Delta' = \lfloor \frac{N}{N'} \rfloor \frac{i}{N}$,
$\widetilde y_j' = b + j \delta'$ and $\widetilde z_k' = b + k \delta'$,
where $m = m_1 m_2$, $m_1 = m_2$, $m' = m_1' m_2'$, $m_1' = m_2'$,
$m = \OO(N)$,  $m' = \OO(N')$, $\delta = \frac{1-2b}{m_1}$,
$\delta' = \frac{1-2b}{m_1'}$,
$\delta' = 2\delta$ and $\Delta' = 4\Delta$.

For the thinned data $\mathbb X_{m,N}^{(b)}$, we define the triple increments
\begin{align*}
T_{i,j,k} X 
&= X_{t_{i}}(\widetilde y_{j},\widetilde z_{k})
-X_{t_{i}}(\widetilde y_{j-1},\widetilde z_{k})
-X_{t_{i}}(\widetilde y_{j},\widetilde z_{k-1})
+X_{t_{i}}(\widetilde y_{j-1},\widetilde z_{k-1})
\\
&\qquad
-X_{t_{i-1}}(\widetilde y_{j},\widetilde z_{k})
+X_{t_{i-1}}(\widetilde y_{j-1},\widetilde z_{k})
+X_{t_{i-1}}(\widetilde y_{j},\widetilde z_{k-1})
-X_{t_{i-1}}(\widetilde y_{j-1},\widetilde z_{k-1}).
\end{align*}
We then obtain for $r \equiv \delta/\sqrt{\Delta} \in (0,\infty)$,
\begin{equation*}
\frac{1}{m N \Delta^\alpha} \sum_{k=1}^{m_2} \sum_{j=1}^{m_1} \sum_{i=1}^N
(T_{i,j,k} X)^2
\pto g_{r,\alpha}(\vartheta),
\end{equation*}
where $g_{r,\alpha}(\vartheta)$ is given by \eqref{eq-pf-1} below.
Let $T_{i,j,k}' X$ be the triple increments for the thinned data 
$\mathbb X_{m',N'}^{(b)}$.
Similarly, we have for $r' \equiv \delta'/\sqrt{\Delta'} \in (0,\infty)$,
\begin{equation*}
\frac{1}{m' N' (\Delta')^\alpha} \sum_{k=1}^{m_2'} \sum_{j=1}^{m_1'} \sum_{i=1}^{N'} 
(T_{i,j,k}' X)^2
\pto g_{r',\alpha}(\vartheta).
\end{equation*}
Since $\delta' = 2 \delta$ and $\Delta' = 4 \Delta$, we obtain $r' = r$ and
\begin{equation*}
\frac{\displaystyle 
\frac{1}{m' N'  (\Delta')^\alpha} \sum_{k=1}^{m_2'} \sum_{j=1}^{m_1'} \sum_{i=1}^{N'} 
(T_{i,j,k}' X)^2}
{\displaystyle
\frac{1}{m N \Delta^\alpha} \sum_{k=1}^{m_2} \sum_{j=1}^{m_1} \sum_{i=1}^N 
(T_{i,j,k} X)^2}
\pto 1.
\end{equation*}
Since it follows that
\begin{align*}
\frac{\displaystyle\frac{1}{m' N' (\Delta')^\alpha} \sum_{k=1}^{m'_2} \sum_{j=1}^{m'_1} 
\sum_{i=1}^{N'} (T_{i,j,k}'X)^2}
{\displaystyle\frac{1}{m N \Delta^\alpha} \sum_{k=1}^{m_2} \sum_{j=1}^{m_1} \sum_{i=1}^N 
(T_{i,j,k}X)^2}
&= \frac{\Delta^\alpha}{(\Delta')^\alpha}
\times
\frac{\displaystyle\frac{1}{m' N'} \sum_{k=1}^{m'_2} \sum_{j=1}^{m'_1} 
\sum_{i=1}^{N'} (T_{i,j,k}'X)^2}
{\displaystyle\frac{1}{m N} \sum_{k=1}^{m_2} \sum_{j=1}^{m_1} \sum_{i=1}^N 
(T_{i,j,k}X)^2}
\\
&= \frac{1}{4^\alpha}
\times
\frac{\displaystyle\frac{1}{m' N'} \sum_{k=1}^{m'_2} \sum_{j=1}^{m'_1} 
\sum_{i=1}^{N'} (T_{i,j,k}'X)^2}
{\displaystyle\frac{1}{m N} \sum_{k=1}^{m_2} \sum_{j=1}^{m_1} \sum_{i=1}^N 
(T_{i,j,k}X)^2},
\end{align*}
we have
\begin{equation*}
\frac{\displaystyle\frac{1}{m' N'} \sum_{k=1}^{m'_2} \sum_{j=1}^{m'_1} 
\sum_{i=1}^{N'} (T_{i,j,k}'X)^2}
{\displaystyle\frac{1}{m N} \sum_{k=1}^{m_2} \sum_{j=1}^{m_1} \sum_{i=1}^N 
(T_{i,j,k}X)^2}
\pto 4^\alpha.
\end{equation*}
We thus define the estimator of the damping parameter $\alpha$ as follows.
\begin{equation*}
\widehat \alpha =
\log \left( \frac{\displaystyle\frac{1}{m' N'} \sum_{k=1}^{m'_2} \sum_{j=1}^{m'_1} 
\sum_{i=1}^{N'} (T_{i,j,k}' X)^2}
{\displaystyle\frac{1}{m N} \sum_{k=1}^{m_2} \sum_{j=1}^{m_1} \sum_{i=1}^N 
(T_{i,j,k}X)^2}
\right)/\log(4).
\end{equation*}

We then obtain the following result.
\begin{thm}\label{th1}
Under [A1], as $m \to \infty$ and $N \to \infty$,
\begin{equation*}
\sqrt{m N} (\widehat \alpha -\alpha^*) = \OO_\PP(1).
\end{equation*}
\end{thm}

\begin{rmk}
\ 
\begin{enumerate}
\item[(1)]
This result is also true for the SPDE driven by 
the $Q_1$-Wiener process in \cite{TKU2023} 
without changing the form of the estimator $\widehat \alpha$.

\item[(2)]
The estimator $\widehat \alpha$ has a faster convergence rate 
than the estimator proposed by \cite{Bossert2023arXiv}.
Substituting this estimator in the contrast function of
\cite{TKU2023} or \cite{TKU2024b}, we will be able to construct a consistent estimator of 
the coefficient parameters of SPDE \eqref{2d_spde} even if $\alpha$ is unknown.
This is future work.
\end{enumerate}
\end{rmk}

\section{Simulations}\label{sec3}
A numerical solution of SPDE \eqref{2d_spde} driven by the $Q$-Wiener process 
\eqref{QW} is generated by
\begin{equation*}
\widetilde X_{t_{i}}(y_j, z_k)
= \sum^{K}_{l_1 = 1}\sum^{L}_{l_2 = 1}x_{l_1,l_2}(t_{i})
e_{l_1}^{(1)}(y_j) e_{l_2}^{(2)}(z_k)
\end{equation*}
with
\begin{equation*}
\dd x_{l_1,l_2}(t) =
-\lambda_{l_1,l_2} x_{l_1,l_2}(t)\dd t
+\sigma \mu_{l_1,l_2}^{-\alpha/2} \dd w_{l_1,l_2}(t),
\quad
x_{l_1,l_2}(0) = \langle \xi, e_{l_1,l_2} \rangle,
\end{equation*}
$i = 1,\ldots, N$, $j = 1,\ldots, M_1$, $k = 1,\ldots, M_2$.
In this simulation, 
we set $(\theta_0, \theta_1,\eta_1, \theta_2, \sigma) = (0,0.2,0.2,0.2,1)$, 
$\alpha^* = 0.5$, $\mu_0 = -19.5$,  $\xi = 0$,
$N = 10^3$, $M_1 = M_2 = 200$, $K = L = 10^4$.
We used R language to compute the estimator $\widehat \alpha$ in Theorem \ref{th1}.
The number of Monte Carlo iterations is 200.

We estimated $\alpha$ using the thinned data 
$\mathbb X_{m,N}^{(b)}$ and $\mathbb X_{m',N'}^{(b)}$
with $b \in \{0.02, 0.04, 0.06, 0.08, 0.1 \}$, 
$m_1 = m_2 \in \{20, 30, \ldots, 90 \}$, $m_1' = m_2' = m_1/2$, $N'=N/4$.
Table \ref{tab1} presents the simulation results of the means 
and the standard deviations of $\widehat \alpha$.
We observe that the biases of $\widehat \alpha$ get smaller 
as $m_1$ and $m_2$ increase regardless of the value of $b$.

\begin{table}[h]
\caption{Means and standard deviations of $\widehat \alpha$. 
The true value $\alpha^* = 0.5$.}
\label{tab1}
\begin{center}
\begin{tabular*}{0.75\textwidth}{@{\extracolsep{\fill}}c|ccccc}
$m_1,m_2$ & & & $b$ & \\ 
 & $0.02$ & $0.04$ & $0.06$ & $0.08$ & $0.1$
\\ \hline
$20$ & $0.430$ & $0.448$ & $0.450$ & $0.451$ & $0.452$
\\
& $(0.01302)$ & $(0.01157)$ & $(0.01222)$ & $(0.01174)$ & $(0.01302)$
\\
$30$ & $0.457$ & $0.468$ & $0.467$ & $0.468$ & $0.471$
\\
& $(0.00853)$ & $(0.00825)$ & $(0.00774)$ & $(0.00794)$ & $(0.00774)$
\\
$40$ & $0.467$ & $0.477$ & $0.481$ & $0.480$ & $0.476$
\\
& $(0.00619)$ & $(0.00592)$ & $(0.00569)$ & $(0.00615)$ & $(0.00619)$
\\
$50$ & $0.476$ & $0.483$ & $0.487$ & $0.489$ & $0.485$
\\
& $(0.00489)$ & $(0.00442)$ & $(0.00497)$ & $(0.00480)$ & $(0.00452)$
\\
$60$ & $0.485$ & $0.485$ & $0.484$ & $0.488$ & $0.492$
\\
& $(0.00365)$ & $(0.00392)$ & $(0.00392)$ & $(0.00365)$ & $(0.00360)$
\\
$70$ & $0.490$ & $0.497$ & $0.501$ & $0.503$ & $0.498$
\\
& $(0.00323)$ & $(0.00321)$ & $(0.00327)$ & $(0.00332)$ & $(0.00334)$
\\
$80$ & $0.502$ & $0.501$ & $0.497$ & $0.492$ & $0.498$
\\
& $(0.00258)$ & $(0.00302)$ & $(0.00270)$ & $(0.00263)$ & $(0.00294)$
\\
$90$ & $0.493$ & $0.492$ & $0.491$ & $0.498$ & $0.506$
\\
& $(0.00251)$ & $(0.00269)$ & $(0.00254)$ & $(0.00236)$ & $(0.00256)$
\end{tabular*}
\end{center}
\end{table}

\section{Proofs}\label{sec4}
We provide the proof of the main result.
\begin{proof}[\bf{Proof of Theorem \ref{th1}}]
For $\delta/\sqrt{\Delta} \equiv r \in(0,\infty)$, we define
\begin{equation*}
\psi_{r,\alpha}(\theta_2)
=\frac{2}{\theta_2^{1-\alpha} \pi}
\int_0^\infty 
\frac{1-\ee^{-x^2}}{x^{1+2\alpha}}
\biggl(
J_0\Bigl(\frac{\sqrt{2}r x}{\sqrt{\theta_2}}\Bigr)
-2J_0\Bigl(\frac{r x}{\sqrt{\theta_2}}\Bigr)+1
\biggr) \dd x,
\end{equation*}
where $J_0$ denotes the Bessel function of the first kind of order $0$.
For $\vartheta = (\kappa, \eta, \theta_2, \sigma^2)$, we define
\begin{equation}\label{eq-pf-1}
g_{r,\alpha}(\vartheta) = 
\frac{\sigma^2 \psi_{r,\alpha}(\theta_2)}{(1-2b)^2} 
\int_{b}^{1-b} \int_{b}^{1-b} \ee^{-\kappa y -\eta z} \dd y \dd z.
\end{equation}
Since 
\begin{equation*}
J_0(\sqrt{2}x)-2J_0(x)+1
= \frac{2}{\pi} \int_0^{\pi/2} (1-\cos(x \cos(t)))(1-\cos(\sin(t))) \dd t \ge 0,
\quad x \ge 0,
\end{equation*}
we note that $g_{r,\alpha}(\vartheta) > 0$.

\textbf{Step 1:}
We first show
\begin{equation}\label{eq-pf-2}
\frac{1}{m N \Delta^{\alpha^*}} \sum_{k=1}^{m_2} \sum_{j=1}^{m_1} \sum_{i=1}^N 
\EE[(T_{i,j,k}X)^2]
= g_{r,\alpha^*}(\vartheta)+ \OO (\Delta).
\end{equation}
We see from (4.20) in \cite{TKU2024b} that
\begin{equation*}
\frac{1}{N \Delta^{\alpha^*}}\sum_{i=1}^N \EE[(T_{i,j,k}X)^2]
= \sigma^2 \ee^{-\kappa \overline y_j -\eta \overline z_k} \psi_{r,\alpha^*}(\theta_2)
+ \OO (\Delta)
\end{equation*}
with $\overline y_j = (\widetilde y_{j-1} +\widetilde y_j)/2$, 
$\overline z_k = (\widetilde z_{k-1} +\widetilde z_k)/2$
for $j= 1, \ldots, m_1$, $k = 1,\ldots, m_2$.
Since it follows that for $f \in C^2((b,1-b))$,
\begin{align*}
&\Biggl| 
\frac{1}{m_1} \sum_{j=1}^{m_1} f(\overline y_j)
- \frac{1}{1-2b} \int_b^{1-b} f(y) \dd y
\Biggr|
\\
&=
\Biggl| \frac{1}{1-2b} 
\sum_{j=1}^{m_1} \int_{\widetilde y_{j-1}}^{\widetilde y_j} 
(f(\overline y_j) -f(y)) \dd y
\Biggr|
\\
&\le 
\Biggl| \frac{1}{1-2b} 
\sum_{j=1}^{m_1} \int_{\widetilde y_{j-1}}^{\widetilde y_j} 
f'(\overline y_j)(y -\overline y_j) \dd y
\Biggr|
\\
&\qquad +
\Biggl| \frac{1}{1-2b} 
\sum_{j=1}^{m_1} \int_{\widetilde y_{j-1}}^{\widetilde y_j} 
\biggl(
\int_0^1 f''(\overline y_j +u(y -\overline y_j)) \dd u 
\biggr)
(y -\overline y_j)^2 \dd y
\Biggr|
\\
&=
0 + \OO \Biggl( \sum_{j=1}^{m_1} 
\int_{\widetilde y_{j-1}}^{\widetilde y_j} (y -\overline y_j)^2 \dd y \Biggr)
\\
&= \OO(\delta^2) = \OO(\Delta),
\end{align*}
we obtain
\begin{align*}
\frac{1}{m N \Delta^{\alpha^*}} \sum_{k=1}^{m_2} \sum_{j=1}^{m_1} \sum_{i=1}^N 
\EE[(T_{i,j,k}X)^2]
&= \sigma^2 \psi_{r,\alpha^*}(\theta_2) \times
\frac{1}{m} \sum_{k=1}^{m_2} \sum_{j=1}^{m_1} 
\ee^{-\kappa \overline y_j -\eta \overline z_k}
+ \OO (\Delta)
\\
&= g_{r,\alpha^*}(\vartheta)+ \OO (\Delta).
\end{align*}

\textbf{Step 2:}
We next show 
\begin{equation}\label{eq-pf-3}
\sqrt{m N} \Biggl(
\frac{1}{m N \Delta^{\alpha^*}} \sum_{k=1}^{m_2} \sum_{j=1}^{m_1} \sum_{i=1}^N 
(T_{i,j,k}X)^2
-g_{r,\alpha^*}(\vartheta)
\Biggr) = \OO_\PP(1).
\end{equation}
We find from Lemma 4.11 in \cite{TKU2024b} that 
\begin{align*}
&\EE \Biggl[ \biggl(
\sum_{k=1}^{m_2} \sum_{j=1}^{m_1} \sum_{i=1}^N 
\Bigl( (T_{i,j,k}X)^2 -\EE[(T_{i,j,k}X)^2] \Bigr)
\biggr)^2 \Biggr]
\\
&= 
\sum_{k,k'=1}^{m_2} \sum_{j,j'=1}^{m_1} \sum_{i,i'=1}^{N}
\cov \Bigl[(T_{i,j,k}X)^2, (T_{i',j',k'}X)^2 \Bigr]
\\
&= \OO(m N \Delta^{2 \alpha^*}),
\end{align*}
which together with \eqref{eq-pf-2} yields
\begin{align*}
&\sqrt{m N} \Biggl(
\frac{1}{m N \Delta^{\alpha^*}} \sum_{k=1}^{m_2} \sum_{j=1}^{m_1} \sum_{i=1}^N 
(T_{i,j,k}X)^2 -g_{r,\alpha^*}(\vartheta)
\Biggr)
\\
&=
\frac{1}{\sqrt{m N} \Delta^{\alpha^*}} \sum_{k=1}^{m_2} \sum_{j=1}^{m_1} \sum_{i=1}^N 
\Bigl( (T_{i,j,k}X)^2 -\EE[(T_{i,j,k}X)^2] \Bigr)
\\
&\qquad +
\sqrt{m N} \Biggl(
\frac{1}{m N \Delta^{\alpha^*}} \sum_{k=1}^{m_2} \sum_{j=1}^{m_1} \sum_{i=1}^N 
\EE[(T_{i,j,k}X)^2] -g_{r,\alpha^*}(\vartheta)
\Biggr)
\\
&=
\OO_\PP(1).
\end{align*}

\textbf{Step 3:}
Finally, we prove
\begin{equation*}
\sqrt{m N}(\widehat \alpha -\alpha^*) = \OO_\PP(1).
\end{equation*}
We define
\begin{equation*}
\mathcal Z =
\frac{1}{m N \Delta^{\alpha^*}} \sum_{k=1}^{m_2} \sum_{j=1}^{m_1} \sum_{i=1}^N 
(T_{i,j,k}X)^2,
\quad
\mathcal Z' =
\frac{1}{m' N' (\Delta')^{\alpha^*}} \sum_{k=1}^{m'_2} \sum_{j=1}^{m'_1} 
\sum_{i=1}^{N'} (T_{i,j,k}' X)^2.
\end{equation*}
In the same way as in Step 2, we have
\begin{equation}\label{eq-pf-4}
\sqrt{m N} \bigl( \mathcal Z' -g_{r,\alpha^*}(\vartheta) \bigr) = \OO_\PP(1).
\end{equation}
Since we have
\begin{equation*}
\biggl| \frac{\mathcal Z'}{\mathcal Z}-1 \biggr|
= \biggl| \frac{\mathcal Z'-g_{r,\alpha^*}(\vartheta) 
-(\mathcal Z -g_{r,\alpha^*}(\vartheta))}
{\mathcal Z -g_{r,\alpha^*}(\vartheta) +g_{r,\alpha^*}(\vartheta)} \biggr|
\le 
\frac{|\mathcal Z'- g_{r,\alpha^*}(\vartheta)| +|\mathcal Z -g_{r,\alpha^*}(\vartheta)|}
{|g_{r,\alpha^*}(\vartheta) -|\mathcal Z -g_{r,\alpha^*}(\vartheta)||},
\end{equation*}
we obtain by \eqref{eq-pf-3} and \eqref{eq-pf-4},
\begin{equation*}
\sqrt{m N} \biggl| \frac{\mathcal Z'}{\mathcal Z}-1 \biggr|
\le 
\frac{\sqrt{m N}|\mathcal Z'- g_{r,\alpha^*}(\vartheta)| 
+\sqrt{m N}|\mathcal Z -g_{r,\alpha^*}(\vartheta)|}
{|g_{r,\alpha^*}(\vartheta) -|\mathcal Z -g_{r,\alpha^*}(\vartheta)||}
= \OO_\PP(1).
\end{equation*}
Using the Taylor expansion
\begin{equation*}
\log(x) = (x-1) - \int_0^1 \frac{1-u}{(1+u(x-1))^2} \dd u (x-1)^2
\end{equation*}
and $\sqrt{m N} (\mathcal Z'/\mathcal Z -1) = \OO_\PP(1)$, we have
\begin{align*}
\sqrt{m N}(\widehat \alpha -\alpha^*)
&= \sqrt{m N} \left(
\log \left( \frac{\displaystyle\frac{1}{m' N'} \sum_{k=1}^{m'_2} \sum_{j=1}^{m'_1} 
\sum_{i=1}^{N'} (T_{i,j,k}'X)^2}
{\displaystyle\frac{1}{m N} \sum_{k=1}^{m_2} \sum_{j=1}^{m_1} \sum_{i=1}^N 
(T_{i,j,k}X)^2}
\right)/\log(4)
-\alpha^*
\right)
\\
&=
\sqrt{m N} \Biggl(
\log \biggl( \frac{(\Delta')^{\alpha^*} \mathcal Z'}
{\Delta^{\alpha^*} \mathcal Z} \biggr)
/\log(4)
-\alpha^*
\Biggr)
\\
&=
\sqrt{m N} \Biggl(
\biggl(
\log \biggl( \frac{\mathcal Z'}{\mathcal Z} \biggr)
+\alpha^* \log(4)
\biggr)/\log(4)
-\alpha^*
\Biggr)
\\
&=
\frac{\sqrt{m N}}{\log(4)}
\log \biggl( \frac{\mathcal Z'}{\mathcal Z} \biggr)
\\
&=
\frac{\sqrt{m N}}{\log(4)}
\Biggl( \frac{\mathcal Z'}{\mathcal Z} -1
+ \OO_\PP \biggl( \frac{1}{m N} \biggr)
\Biggr)
\\
&=\OO_\PP(1).
\end{align*}
This concludes the proof.
\end{proof}


\begin{thebibliography}{99}



\bibitem{Bibinger_Bossert2023}
M. Bibinger and P. Bossert. 
\newblock Efficient parameter estimation for parabolic SPDEs based on 
a log-linear model for realized volatilities. 
\newblock {\em Japanese Journal of Statistics and Data Science}, 6:407--429, 2023.


\bibitem{Bibinger_Trabs2020}
M.~Bibinger and M.~Trabs.
\newblock Volatility estimation for stochastic {PDE}s using high-frequency
  observations.
\newblock {\em Stochastic Processes and their Applications}, 130(5):3005--3052, 2020.


\bibitem{Bossert2023arXiv}
P.~Bossert.
\newblock Parameter estimation for second-order SPDEs in multiple space dimensions.
\newblock{ \em arXiv preprint arXiv:2310.17828}, 2023.


\bibitem{Chong2020}
C. Chong. 
\newblock High-frequency analysis of parabolic stochastic pdes. 
\newblock{ \em The Annals of Statistics}, 48(2):1143--1167, 2020.


\bibitem{Cialenco_Huang2020}
 I. Cialenco and Y. Huang. 
\newblock A note on parameter estimation for discretely sampled spdes. 
\newblock{ \em Stochastics and Dynamics}, 20(3):2050016. 2020.


\bibitem{Gaudlitz_Reiss2023}
S. Gaudlitz and M. Reiss. 
\newblock Estimation for the reaction term in semi-linear SPDEs
under small diffusivity.
\newblock {\em Bernoulli} 29(4):3033-3058, 2023.


\bibitem{Hildebrandt_Trabs2021}
F.~Hildebrandt and M.~Trabs.
\newblock Parameter estimation for {SPDE}s based on discrete observations in
  time and space.
\newblock {\em Electronic Journal of Statistics}, 15(1):2716--2776, 2021.


\bibitem{Kaino_Uchida2020}
Y.~Kaino and M.~Uchida.
\newblock Parametric estimation for a parabolic linear {SPDE} model based on
  discrete observations.
\newblock {\em Journal of Statistical Planning and Inference}, 211:190--220, 2020.


\bibitem{Kaino_Uchida2021}
Y.~Kaino and M.~Uchida.
\newblock Adaptive estimator for a parabolic linear {SPDE} with a small noise.
\newblock {\em Japanese Journal of Statistics and Data Science}, 4:513--541, 2021.


\bibitem{TKU2023}
Y.~Tonaki, Y.~Kaino, and M.~Uchida.
\newblock Parameter estimation for linear parabolic {SPDEs} in two space
  dimensions based on high frequency data.
\newblock {\em Scandinavian Journal of Statistics}, 50(4):1568-1589, 2023.


\bibitem{TKU2024a}
Y.~Tonaki, Y.~Kaino, and M.~Uchida.
\newblock Parameter estimation for a linear parabolic {SPDE} model in two space
  dimensions with a small noise.
\newblock {\em Statistical Inference for Stochastic Processes}, 27(1):123--179, 2024.


\bibitem{TKU2024b}
Y.~Tonaki, Y.~Kaino, and M.~Uchida.
\newblock Parametric estimation for linear parabolic {SPDE}s in two space
  dimensions based on temporal and spatial increments.
\newblock {\em Metrika}, 2024. \url{https://doi.org/10.1007/s00184-024-00969-x}


\bibitem{TKU2024arXiv}
Y.~Tonaki, Y.~Kaino, and M.~Uchida.
\newblock Small diffusivity asymptotics for a linear parabolic {SPDE} 
  in two space dimensions.
\newblock{ \em arXiv preprint arXiv:2404.02513}, 2024.
\end{thebibliography}
\end{document}